\documentclass[12pt]{article}

\usepackage{amsmath,amsthm,amsfonts,amssymb,amscd,epsf,epsfig,psfrag,enumerate}
\usepackage[mathscr]{eucal}

\title{Mapping Class Groups of Compression Bodies \break and 3-Manifolds}

\author{ Ulrich Oertel}

\date{January, 2006; revised June, 2007}

\newtheorem{thm}{Theorem}[section] \newtheorem{lemma}[thm]{Lemma}
 
\newtheorem{claim}[thm]{Claim}
\newtheorem{proposition}[thm]{Proposition}
 \newtheorem*{claim*}{Claim}

 \theoremstyle{definition}
\newtheorem{defn}[thm]{Definition}
\newtheorem{defns}[thm]{Definitions} 
 \newtheorem{ex}[thm]{Example}
 \newtheorem{remark}[thm]{Remark}
\theoremstyle{remark}

\def\figurefontsize{12}

\begin{document}

\maketitle

\vskip 0.5in

\def\HDS{half-disc sum}

\def\Im{\text{Im}}
\def\im{\text{Im}}
\def\rel{\text{ rel }}
\def\irred{irreducible}
\def\half{spinal pair }
\def\spinal{\half}
\def\spinals{\halfs}
\def\halfs{spinal pairs }
\def\reals{\mathbb R}
\def\rationals{\mathbb Q}
\def\complex{\mathbb C}
\def\naturals{\mathbb N}
\def\integers{\mathbb Z}
\def\id{\text{id}}

\def\proj{P}
\def\hyp {\hbox {\rm {H \kern -2.8ex I}\kern 1.15ex}}

\def\Diff{\text{Diff}}

\def\weight#1#2#3{{#1}\raise2.5pt\hbox{$\centerdot$}\left({#2},{#3}\right)}
\def\intr{{\rm int}}
\def\inter{\ \raise4pt\hbox{$^\circ$}\kern -1.6ex}
\def\Cal{\cal}
\def\from{:}
\def\inverse{^{-1}}
\def\Max{{\rm Max}}
\def\Min{{\rm Min}}
\def\fr{{\rm fr}}
\def\embed{\hookrightarrow}
\def\Genus{{\rm Genus}}
\def\Z{Z}
\def\X{X}
\def\interior{\text{int}}

\def\roster{\begin{enumerate}}
\def\endroster{\end{enumerate}}
\def\intersect{\cap}
\def\definition{\begin{defn}}
\def\enddefinition{\end{defn}}
\def\subhead{\subsection\{}
\def\theorem{thm}
\def\endsubhead{\}}
\def\head{\section\{}
\def\endhead{\}}
\def\example{\begin{ex}}
\def\endexample{\end{ex}}
\def\ves{\vs}
\def\mZ{{\mathbb Z}}
\def\M{M(\Phi)}
\def\bdry{\partial}
\def\hop{\vskip 0.15in}
\def\jerk{\vskip 0.08in}
\def\mathring{\inter}
\def\trip{\vskip 0.09in}
\def\H{\mathscr{H}}
\def\S{\mathscr{S}}
\def\E{\mathscr{E}}
\def\K{\mathscr{K}}
\def\B{\mathscr{B}}
\def\cl{\text{Cl}}

\begin{abstract} We analyze the mapping class group $\H_x(W)$ of automorphisms of the exterior boundary $W$ of a compression
body $(Q,V)$ of dimension 3 or 4 which extend over the compression body.  Here $V$ is the interior boundary of the compression body $(Q,V)$.  Those
automorphisms which extend as automorphisms of
$(Q,V)$ rel
$V$ are called discrepant automorphisms, forming the mapping class group
$\H_d(W)$ of discrepant automorphisms of $W$ in $Q$.  If $\H(V)$ denotes the mapping class group of $V$ we describe a short exact sequence of
mapping class groups relating $\H_d(W)$, $\H_x(W)$, and $\H(V)$.

For an orientable, compact, reducible 3-manifold $W$, there is a canonical ``maximal" 4-dimensional compression body
whose exterior boundary is
$W$ and whose interior boundary is the disjoint union of the irreducible summands of $W$.  Using the canonical compression body
$Q$, we show that the mapping class group
$\H(W)$ of
$W$ can be identified with
$\H_x(W)$.  Thus we obtain a short exact sequence for the mapping class group of a 3-manifold, which gives
the mapping class group of the disjoint union of irreducible summands as a quotient of the entire mapping class
group by the group of discrepant automorphisms.  The group of discrepant automorphisms is described in terms of
generators, namely certain ``slide" automorphisms, ``spins," and ``Dehn twists" on 2-spheres.

Much of the motivation for the results in this paper comes from a research
program for classifying automorphisms of compact 3-manifolds, in the spirit of the Nielsen-Thurston classification of automorphisms of
surfaces.  
\end{abstract}

\section{Introduction}\label{Intro}

\begin{figure}[ht]
\centering
 \psfrag{S}{\fontsize{\figurefontsize}{12}$S$}\psfrag{H}{\fontsize{\figurefontsize}{12}$H$}\psfrag{a}{\fontsize{\figurefontsize}{12}$\alpha$}

\scalebox{1.0}{\includegraphics{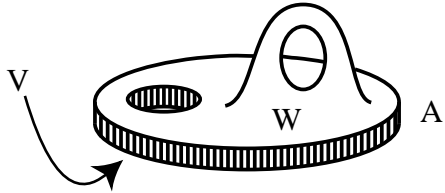}} \caption{\small
Compression body.} \label{MapCompression}
\end{figure}

\begin{definition}  A {\it compression body} is a manifold triple $(Q,W,V)$ of any dimension $n\ge 3$ constructed as
$V\times I$ with 1-handles attached to $V\times 1$, and with $V=V\times 0$.   Here $V$ is a manifold, possibly with
boundary, of dimension $n-1\ge 2$.  We use the symbol $A$ to denote $\bdry
V\times I$ and
$W$ to denote
$\bdry Q-(V\cup
\intr(A))$.  See Figure \ref{MapCompression}.  When $W$ is closed, we can regard the compression body as a manifold pair, and even
when $\bdry W\ne\emptyset$, we will often denote the compression body simply as $(Q,V)$.  
\end{definition}

\noindent{\bf Remark:}  In most applications, a compression body satisfies the condition that $V$ contains no sphere components, but we choose not to
make this a part of the definition.

\hop

We note that, when $\bdry V\ne\emptyset$, it is crucial in this paper that there be a product ``buffer zone" $A$ between $\bdry V$ and
$\bdry W$.  Some applications of the isotopy extension theorem would fail without this buffer. 

\begin{definition}  The symbol $\H$ will be used to denote the mapping class group of orientation preserving
automorphisms of an orientable manifold, a manifold pair, or a manifold triple.  When $(Q,V)$ is a compression body, and $\bdry
V=\emptyset$,
$\H(Q,V)$ is the mapping class group of the manifold pair; when $\bdry V\ne\emptyset$ by abuse of notation, we use $\H(Q,V)$ to
denote $\H(Q, W,V)$.
\end{definition}

\noindent{\bf Convention:} Throughout this paper, manifolds will be compact and orientable; automorphisms will be orientation preserving.

\begin{definition}\label{BasicDefinitions} Suppose $(Q,V)$ is a compression body.  Let
$\H_x (W)$ denote the image of the restriction homomorphism $H(Q,V)\to \H(W)$, i.e. the subgroup of $\H(W)$ consisting of automorphisms of
$W$ which extend over
$Q$ to yield an element of $\H (Q,V)$.  We call the elements of $\H_x (W)$ {\it extendible} automorphisms of $W$.  
Let $\H_d(W)$ denote the image of of the restriction map $H(Q 
\rel V)\to \H(W)$, i.e. the subgroup of $\H_x(W)$ consisting of
isotopy classes of automorphisms
$f$ of
$W$ which extend over $Q$ such that $f|_V$ becomes the identity.  We call the elements of $\H_d (W)$ {\it discrepant} automorphisms of $W$.  
\end{definition}

A fundamental fact well known for compression bodies of dimension 3 is that an automorphism of $W$ which
extends over a compression body $Q$ with $\bdry_eQ=W$ uniquely determines up to isotopy the restriction of the
extended automorphism to $V=\bdry_iQ$.  In general we have a weaker statement:

\begin{thm}\label{FundamentalTheorem} Suppose $(Q,V)$ is a compression body of dimension $n\ge 3$.  Suppose at most one component of $V$ is
homeomorphic to a sphere
$S^{n-1}$ and suppose each component of
$V$ has universal cover either contractible or homeomorphic to $S^{n-1}$.
If $f:(Q,V)\to (Q,V)$ is an automorphism with
$f|_W$ isotopic to the identity, then:
\jerk

\noindent (a) {$f|_V$ is homotopic to the identity.  It follows that for an automorphism $f$ of $(Q,V)$, 
$f|_W$ determines
$f|_V$ up to homotopy.}
\jerk

\noindent (b) {If $Q$ is 3-dimensional, $f|_W$ determines $f|_V$ up to isotopy, and $\H_x(W)\approx \H(Q,V)$. }
\jerk

\noindent (c) {If $Q$ is 4-dimensional, and each component of $V$ is irreducible, then $f|_W$ determines $f|_V$ up to isotopy.}
\end{thm}

Suppose $(Q,V)$ is a compression body of dimension 3 or 4  satisfying the conditions of Theorem \ref{FundamentalTheorem} (b) or (c).  Then
there is an {\it eduction} homomorphism $e: \H_x(W)\to \H(V)$.  For $f\in  \H_x(W)$, we extend $f$ to $\bar
f:(Q,V)\to (Q,V)$, then we restrict $\bar f$ to $V$ obtaining $e(f):V\to V$.  By the theorem, $e$ is well-defined.

\begin{thm}\label{SequenceThm} Suppose $(Q,V)$ is a compression body of dimension 3 or 4.  Suppose at most one component of $V$ is a sphere. In case
$Q$ has dimension 4, suppose every component of 
$V$ is a 3-manifold whose universal cover is contractible or $S^3$. Then the 
following is a short exact sequence of mapping class groups.
$$1\to \H_d(W)\to \H_x(W)\xrightarrow{e} \H(V)\to 1.$$

\noindent The first map is inclusion, the second map is eduction.
\end{thm}

There is a {\it canonical 4-dimensional compression body} $(Q,V)$ associated to a compact reducible 3-manifold $W$ without
boundary spheres, see Section \ref{ThreeManifoldSection}.  It has the following properties.

\begin{proposition}\label{UniquenessProposition}  The canonical compression body $(Q,V)$ associated to a compact reducible 3-manifold $W$ has exterior
boundary
$W$ and interior boundary the disjoint union of its irreducible summands.  It is unique in the following sense:  If $(Q_1,W,V_1)$ and $(Q_2,W,V_2)$ are
two canonical compression bodies, and $v:V_1\to V_2$ is any homeomorphism, then there is a homeomorphism $g:(Q_1,V_1)\to (Q_2,V_2)$ 
with $g|_{V_1}=v$.
\end{proposition}

We will initially construct a canonical compression body using a ``marking" by a certain type of system $\S$ of essential spheres in
$W$, called a {\it symmetric system}.   The result is the
{\it canonical compression body associated to $W$ with the symmetric system $\S$}.  The canonical compression body associated to $W$ with $\S$ is
unique in a different sense: If $(Q_1,W,V_1)$ and $(Q_2,W,V_2)$ are canonical compression bodies associated to $W$ with $\S$, then the identity on
$W$ extends to a homeomorphism $g:Q_1\to Q_2$, and $g|_{V_1}$ is unique up to isotopy.  Only later will we recognize the uniqueness described in
Proposition \ref{UniquenessProposition} when we compare canonical compression bodies constructed from different symmetric systems of spheres.

There are well-known automorphisms of
$W$ which, as we shall see, lie in
$\H_d(W)$ with respect to the canonical compression body $Q$.  The most important ones are called {\it slides}.  Briefly, a slide
is an automorphism of
$W$ obtained by cutting on an essential sphere $S$ to obtain $W|S$ ($W$ cut on $S$), capping one of the duplicate
boundary spheres coming from $S$ with a ball $B$ to obtain $W'$, sliding the ball along a closed curve in this capped manifold and
extending the isotopy of $B$ to the manifold, then removing
$B$ and reglueing on the two boundary spheres.  An {\it interchanging slide} is done using
two separating essential spheres cutting homeomorphic manifolds from $W$.  Cut on the two spheres $S_1$ and $S_2$, remove the
homeomorphic submanifolds, cap the two boundary spheres with two balls $B_1$ and $B_2$, choose two paths to move $B_1$ to $B_2$
and $B_2$ to $B_1$, extend the isotopy, remove the balls, and reglue.  Interchanging slides are used to interchange
prime summands.   A {\it spin} does a  ``half rotation" on an
$S^2\times S^1$ summand.   In addition, for each embedded essential sphere, there is an automorphism of at most order 2 called a {\it Dehn twist}
on the  sphere.  If slides, Dehn twists, interchanging slides, and spins are done using spheres from a given symmetric system $\S$ (or
certain separating spheres associated to non-separating spheres of $\S$) then they are called {\it
$\S$-slides},
{\it $\S$-Dehn twists},
{\it $\S$-interchanging slides}, and
{\it $\S$-spins}.  An $\S$-spin is a ``half-Dehn-twist" on a separating sphere associated to a non-separating
sphere of $\S$, not on a sphere of $\S$.

\begin{proposition}\label{CharacterizationProp}  If $Q$ is the canonical compression body associated to the compact
orientable, 3-manifold
$W$ with symmetric system $\S$ and $W=\bdry_eQ$, $V=\bdry_iQ$, then  
\jerk

\noindent (a) $\H_x(W)=\H(W)$.
\jerk

\noindent (b)  The mapping class group $\H_d(W)$ is generated by $\S$-slides, $\S$-Dehn twists, $\S$-slide
interchanges of $S^2\times S^1$ summands, and $\S$-spins.
\end{proposition}

Proposition \ref{CharacterizationProp} (b) is a version of a standard result due to E. C{\'e}sar de S{\'a}
\cite{EC:Automorphisms}, see also M. Scharlemann in Appendix A of \cite{FB:CompressionBody}, and
D. McCullough in \cite{DM:MappingSurvey} (p. 69). 

From Proposition \ref{CharacterizationProp} and Theorem \ref{SequenceThm} we obtain:

\begin{thm} \label{ThreeManifoldSequenceThm}  If $W$ is a compact, orientable, reducible 3-manifold, and $V$ is the disjoint
union of its irreducible summands, then the following sequence is exact:

$$1\to \H_d(W)\to \H(W)\xrightarrow{e} \H(V)\to 1.$$

\noindent Here $\H_d(W)$ is the group of discrepant automorphisms for a canonical compression body, the first map is inclusion, and the second map
is eduction with respect to a canonical compression body.
\end{thm}

Using automorphisms of compression bodies, one can reinterpret some other useful results.  Suppose
$W$ is an
$m$-manifold ($m=2$ or $m=3$) with some sphere boundary components.  Let
$V_0$ denote the manifold obtained by capping the sphere boundary components with balls and let $P$ denote the union of capping
balls in $V_0$.  There is another appropriately chosen canonical compression body $(Q,V)$ for $W$, such that $V$ is the disjoint
union $V_0\sqcup P$. With respect to this compression body we have a group of discrepant automorphisms, $\H_d(W)$.  We say the
pair $(V_0, P)$ is a {\it spotted manifold}.  The mapping class group of this spotted manifold, as a pair is the same as
$\H(W)$.  The following result probably exists in the literature in some form, and can be proven more directly, and with other terminology.  
For example, a special case is mentioned
in
\cite{JB:Braids}.  The reason for including it here is to make an analogy with Theorem \ref{SequenceThm}.  Let
$f:W\to W$ be an automorphism.  Then $f$ clearly induces a {\it capped automorphism}
$f_c:V_0\to V_0$.   Simply extend $f$ to the balls of $P$ to obtain $f_c$.  Then $f_c$ is determined up to isotopy by $f$.

\begin{thm}\label{SpottedThm}  With $W$ and $(Q,V)$ as above, and $V_0$ a connected surface or a connected irreducible 3-manifold,
there is an exact sequence 

$$1\to \H_d(W)\to \H(W)\xrightarrow{e} \H(V_0)\times \H(P)\to 1.$$

\noindent $\H(P)$ is the group of
permutations of the balls of $P$; $\H(V_0)$ is the mapping class group of the capped manifold.  The group $\H_d(W)$ of discrepant 
automorphisms is the subgroup of automorphisms $f$ of $W$ which map each boundary sphere to itself such that the capped automorphism $f_c$ is isotopic
to the identity.  The map
$\H_d(W)\to
\H(W)$ is inclusion.  The eduction map $e:\H(W)\to \H(V_0)\times \H(P)$ takes $f\in \H(W)$ to $(f_c,\rho)$ where $\rho$ is the
permutation of
$P$ induced by $f$.
\end{thm}

In Section \ref{Sequence} we prove a lemma about mapping class groups of triples.  This is the basis for the results about mapping class
groups of compression bodies. Let
$M$ be an
$n$-ma\-ni\-fold, and  let $W$ and $V$ be disjoint compact submanifolds of any (possibly different) dimensions $\le n$.
The submanifolds $W$ and
$V$ need not be properly embedded (in the sense that $\bdry W\embed \bdry M$) but they must be submanifolds for
which the isotopy extension theorem applies, for example smooth submanifolds with respect to some smooth structure.
In this more general setting, we make some of the same definitions as in Definition \ref{BasicDefinitions}, as well as some new ones.

\begin{defns}  Let $(M,W,V)$ be a triple as above.
We define $\H_x(V)$, the group of {\it extendible automorphisms} of $V$, as the image of the restriction map $\H(M,W,V)\to \H(V)$.   Similarly, we
define
$\H_x(W,V)$, the group of {\it extendible automorphisms} of the pair $(W,V)$, as the image of the restriction map $\H(M,W,V)\to \H(W,V)$, and the group
$\H_x (W)$ of {\it extendible automorphisms} of $W$ as the image of the restriction homomorphism $\H(M,W,V)\to \H(W)$.   The group $\H_d(W)$ of
{\it discrepant} automorphisms of $W$ is the image of the restriction map $\H((M,W) \rel V)\to \H(W)$. 
\end{defns}

\begin{lemma}\label{SequenceGeneralLemma} Suppose $(M,W,V)$ is a triple as above.  The
following is a short exact sequence of mapping class groups.
$$1\to \H_d(W)\to \H_x(W,V)\to \H_x(V)\to 1.$$

\noindent The first map is inclusion of an element of $\H_d(W)$ extended to $(W,V)$ by the identity on $V$, the second map is restriction.
\end{lemma}

 Let $M$ be an irreducible,
$\bdry$-reducible 3-manifold.  There is a 3-dimensional characteristic compression body embedded in $M$ and containing $\bdry M$ as its exterior
boundary, first described by Francis Bonahon, which is unique up to isotopy, see
\cite{FB:CompressionBody}, or see \cite{UO:Autos}, \cite{OC:Autos}.  Removing this from $M$ yields an irreducible
$\bdry$-irreducible 3-manifold.  The following is an application of Lemma \ref{SequenceGeneralLemma} which describes how mapping
class groups interact which the characteristic compression body.

\begin{thm}\label{MappingCharacteristicThm} Let $M$ be an irreducible, $\bdry$-reducible 3-manifold with characteristic compression
body
$Q$ having exterior boundary $W=\bdry M$ and interior boundary $V$.  Cutting $M$ on $V$ yields $Q$ and an irreducible, $\bdry$-irreducible
manifold
$\hat M$.  There is an exact sequence

$$1\to \H_d(W)\to \H_x(W)\to \H_x(V)\to 1,$$

\noindent where $ \H_d(W)$ denotes the usual group of discrepant automorphisms for $Q$, $\H_x(W)$ denotes the group of automorphisms
which extend to $M$ and $\H_x(V)$ denotes the automorphisms of $V$ which extend to $\hat M$.  
\end{thm}

To a large degree, this paper was motivated by a research program of the author and Leonardo N. Carvalho, to classify
automorphisms of compact 3-manifolds in the spirit of the Nielsen-Thurston classification.  This began with a paper of the author
classifying automorphisms of 3-dimensional handlebodies and compression bodies, see \cite{UO:Autos}, and continued with a paper of
Carvalho, \cite{LNC:Tightness}.  Recently, Carvalho and the author completed a paper giving a classification of
automorphisms for compact 3-manifolds, see \cite {OC:Autos}.  At this writing, the paper \cite {OC:Autos} needs to be revised to make it
consistent with this paper.  In the case of a reducible manifold
$W$, the classification applies to elements of $\H_d(W)$ and $\H(V)$, where $V$ is the disjoint union of the irreducible summands of $W$, as
above.  The elements of $\H_d(W)$ are represented by automorphisms of $W$ each having support in a connected sum embedded in $W$ of $S^2\times S^1$'s
and handlebodies.  It is these automorphisms, above all, which are remarkably difficult to understand from the point of view of finding dynamically
nice representatives.  The theory developed so far describes these in terms of automorphisms of compression bodies and handlebodies. Throughout the
theory, both in \cite{UO:Autos} and in \cite {OC:Autos}, the notion of a ``spotted manifold" arises as an irritating technical detail, which we point
out, but do not fully explain.  Theorem \ref{SpottedThm} gives some explanation of the phenomenon.

I thank Leonardo N. Carvalho and Allen Hatcher for helpful discussions related to this paper.

\section{Proof of Theorem \ref{FundamentalTheorem}.}\label{Fundamental}

\begin{proof}[Proof of Theorem \ref{FundamentalTheorem}]  Suppose $(Q,V)$ is a compression body of dimension $n\ge 3$.  We suppose that
every  component $V_i$ of $V$ has universal cover either contractible or homeomorphic to $S^{n-1}$.  We suppose also that at most one component of
$V$ is homeomorphic to a sphere.    A compression body structure for $Q$ is defined by the system $\E$ of $(n-1)$-balls together with a product structure on some
components of $Q|\E$, and an identification of each of the remaining components of $Q|\E$ with the $n$-ball.  As elsewhere in this paper, we choose a
particular type of compression body structure corresponding to a system of balls $\E$ such that exactly one duplicate $E_i'$ of a ball $E_i$ of $\E$
appears on
$V_i\times 1$ for each product component of $V_i\times I$ of $Q|\E$.  In terms of handles, this means we are regarding
$Q$ as being obtained from $n$-balls and products by attaching just one 1-handle to each $V_i\times I$.  For each
component $V_i$ of $V$ we have an embedding $V_i\times I\embed Q$.

Suppose $f:(Q,V)\to (Q,V)$ is the identity on $W$. We will eventually show that for every component $V_0$ of $V$, $f|_{V_0}$ is homotopic to the
identity on $V_0$ (via a homotopy of pairs preserving the boundary). If $V$ contains a single sphere component and
$V_0$ is this sphere component, then clearly $f(V_0)=V_0$ and $f|_{V_0}$ is homotopic to the identity, since every degree 1 map of $S^{n-1}$ is
homotopic to the identity.   

If $V_0$ is not a sphere, we will begin by showing that $f(V_0)=V_0$.  Consider the case where $\pi_1(V_0)$ is non-trivial.  Clearly $\pi_1(Q)$ can be
expressed as a free product of the 
$\pi_1(V_i)$'s and some $\integers$'s.  If $f$ non-trivially permuted $V_i$'s homeomorphic to $V_0$, the induced map $f_*:\pi_1(Q)\to \pi_1(Q)$ would
not be the identity, since it would conjugate factors of the free product.  On the other hand $f_*$ must be the identity, since the inclusion $j:W\to
Q$ induces a surjection on $\pi_1$ and $f:W\to W$ induces the identity on $\pi_1$.  Next we consider the case that
$V_0$ is contractible.  Then $\bdry V_0\ne \emptyset$.  In this case also, we see that $f(V_0)=V_0$, since $f$ is the identity on
$(\bdry V_0\times 1)\subset \bdry W$, and so $f(\bdry V_0)=\bdry V_0$.  From our assumptions, it is always the case that $V_0$ is contractible, is a
sphere, or has non-trivial $\pi_1$.  We conclude that 
$f(V_0)=V_0$ in all cases.

It remains to show that if $V_0$ is not a sphere, then $f|_{V_0}$ is homotopic to the identity.  When $Q$ has dimension 3, this is a standard
result, see for example \cite{UO:Autos}, \cite{KJ:Combinatorics3Mfds}.  The method of proof, using ``vertical essential surfaces" in $V_0\times
I$, depends on the fact that $V_0$ is not a sphere.  In fact, the method shows that  $f|_{V_0}$ is isotopic to the identity, which proves the theorem in
this dimension.

Finally now it remains to show that $f|V_0$ is homotopic to the identity when $n\ge 4$, and $V_0$ is not a sphere.   To do this, we will work
in a cover of
$Q$.  We may assume $V_0$ is either finitely covered by the sphere, but not the sphere, or it has a contractible universal cover.  We will denote by
$\tilde Q_0$ the cover of $Q$ corresponding to the subgroup $\pi_1(V_0)$ in $\pi_1(Q)$.  The compression body $Q$ can be constructed from a disjoint
union of connected product components and balls, joined by 1-handles.  Covering space theory then shows that the cover
$\tilde Q_0$ can also be expressed as a disjoint union of connected products and balls, joined by 1-handles.   The balls and 1-handles are lifts of the
balls and 1-handles used to construct
$Q$.  The products involved are:  one product of the form $V_0\times I$, products of the form
$S^{n-1}\times I$, and products of the form
$U\times I$, where
$U$ is contractible, these being connected components of the preimages of the connected products used to construct $Q$.  There is a distinguished
component of the preimage of
$V_0$ in
$\tilde Q_0$ which is homeomorphic to $V_0$; we denote it $V_0$ rather than $\tilde V_0$, since the latter would
probably be interpreted as the entire preimage of $V_0$ in $\tilde Q_0$.  If $V_0$ has non-trivial fundamental group, then this is
the only component of the preimage of $V_0$ homeomorphic to $V_0$, otherwise it is distinguished simply by the construction of
the cover corresponding to $\pi_1(V_0)$, via choices of base points.  Now let
$\breve Q_0$ denote the manifold obtained from $\tilde Q_0$ by capping all $S^{n-1}$ boundary spheres by $n$-balls. 
Evidently, $\breve Q_0$ has the homotopy type of $V_0$.  This can be seen, for example, by noting that $\tilde Q_0$ is
obtained from $V_0\times I$ and a collection of $n$-balls by using 1-handles to attach universal covers of $V_i\times
I$, $i\ne 0$.  The universal covers of $V_i\times I$ have the form $S^{n-1}\times I$ or $U\times I$ where $U$ is
contractible.  When we cap the boundary spheres of $\tilde Q_0$ to obtain $\breve Q_0$, the $(S^{n-1}\times I)$'s are
capped and become balls.  

Having described the cover $\tilde Q_0$, and the capped cover $\breve Q_0$, we begin the argument for showing that if
$f:(Q,V)\to (Q,V)$ is an automorphism which extends the identity on $W$, and $V_0$ is a non-sphere component of $V$, then $f|_{V_0}$ is homotopic to the
identity.  We have chosen a compression body structure for $Q$
in terms of $\E$, the product $V\times I$, and the collection of balls of $Q|\E$. Using the automorphism $f$
we obtain another compression body structure:  $Q|f(\E)$ contains the product $f(V\times I)$, which includes
$f(V_0\times I)$ as a component.  Recall that we chose $\E$
such that exactly one duplicate $E_0'$ of exactly one ball $E_0$ of
$\E$ appears on $V_0\times I$.  We can lift the inclusion map
$i:V_0\times I\to Q$ to an embedding
$\tilde i: V_0\times I\embed
\tilde Q_0\subset \breve Q_0$.  Again we abuse notation by denoting the image product as
$V_0\times I$, and we use $E_0$ to denote the lift of $E_0$ as well.  Similarly we lift  $f\circ i:V_0\times I\to
Q$ to an embedding
$\widetilde{f\circ i}:V_0\times I\to
\tilde Q_0\subset \breve Q_0$. 

\begin{claim} The lifts $\tilde i$ and $\widetilde
{f\circ i}$ coincide on $(V_0\times 1)-\intr(E_0)$.  In other words, $\tilde f$ is the identity on $(V_0\times 1)-\intr(E_0)$.
\end{claim}

\begin{proof}[Proof of Claim]  We consider two cases: 

\noindent{\it Case 1:   $\bdry V_0\ne \emptyset$.}    The lift $\tilde f$ of $f$ must take $\bdry (V_0\times I)-\intr(E_0)$ in $\tilde Q_0$
to itself by  uniqueness of lifting, since in this case  $\bdry (V_0\times I)-\intr(E_0)$ is connected.    

\noindent{\it Case 2:  $V_0$ has nontrivial fundamental group.}  In this case, lifting $f\circ i$ and $i$, we see that $\widetilde {f\circ
i}(V_0\times 1-E_0)$ and $\tilde {i}(V_0\times 1-E_0)$ are components in the cover $\tilde
Q_0$ of the preimage of $V_0\times 1-E_0\subset W\subset Q$.  From the description
of the cover $\tilde Q_0$, there is only one such component homeomorphic to $V_0-E_0$, hence both lifts must map
$(V_0\times 1)-E_0$ to the same component, and the lifts $\widetilde{f\circ i}$ and $\tilde i$ must coincide on $(V_0\times 1)-\intr(E_0)$ by
uniqueness of lifting.

These two cases cover all the possibilities for $V_0$, for if $V_0$ has trivial fundamental group, then $V_0$ has a contractible
cover by assumption, since we are ruling out the case $V_0$ a sphere.   But this implies the boundary is non-empty, and we apply Case 1. 
Otherwise, we apply Case 2.  

This completes the proof of the claim.
\end{proof}

Following our convention of writing $V_0\times I$ for $\tilde i(V_0\times I)$, we now conclude that though in general $\tilde i|_{E_0}\ne
\widetilde {f\circ i}|_{E_0}$, these lifts are equal on $\bdry E_0$.

\begin{claim} The map $\widetilde{f\circ i}|_{E_0}$ is homotopic rel boundary to $\tilde i|_{E_0}$ in $\breve Q_0$.
\end{claim}

\begin{proof}[Proof of Claim]  We consider two cases; in Case 1 the universal cover of $V_0$ is contractible and in
Case 2 the universal cover of $V_0$ is $S^{n-1}$.  Let $C\embed \tilde W$ denote $\bdry E_0$, an $(n-2)$-sphere, and
$\tilde W_0$ denote the preimage of $W$ in $\tilde Q_0\subset \breve Q_0$, and choose a base point $w_0$ in $C$.  

\hop
\noindent{\it Case 1: The universal cover of $V_0$ is contractible.}  In this case $\breve Q_0$ deformation retracts to $V_0\subset V_0\times I$ and
the universal cover $\widetilde{\breve Q}_0$ of $\breve Q_0$ is contractible, by the discussion above.  Hence all higher homotopy groups of $\breve Q_0$
are trivial.  Using a trivial application of the long exact sequence for relative homotopy groups we obtain:

$$\cdots \pi_{n-1}(\breve Q_0,w_0)\to \pi_{n-1}(\breve Q_0,C,w_0)\to \pi_{n-2}(C,w_0)\to \pi_{n-2}(\breve
Q_0,w_0)\cdots  $$

or

$$\cdots 0\to \pi_{n-1}(\breve Q_0,C,w_0)\to\integers\to 0 \cdots,  $$

\noindent which implies that all maps of $(n-1)$-balls to $(\breve Q_0,C)$ which are degree one on the boundary are
homotopic.

\hop

\noindent{\it Case 2: The universal cover of $V_0$ is $S^{n-1}$ (but $V_0$ is not homeomorphic to $S^{n-1}$).}  In this case again $\breve
Q_0$ deformation retracts to
$V_0\subset V_0\times I\subset \breve Q_0$ by the discussion above.  The universal cover $\widetilde{\breve Q}_0$ of $\breve Q_0$ is homotopy equivalent
to $S^{n-1}$, and it has the structure of $S^{n-1}\times I$ attached by 1-handles to products of the form
$U\times I$ with $U$ contractible, and $n$-balls, some of the $n$-balls coming from capped preimages of $V_i\times I$
where $V_i$ is covered by $S^{n-1}$.  Thus the higher homotopy groups of $\breve Q_0$ are the same as those of $\widetilde{\breve Q}_0$ or $S^{n-1}$. 
Again using an application of the long exact sequence for relative homotopy groups we obtain:

$$\cdots\pi_{n-1}(C,w_0)\to \pi_{n-1}(\breve Q_0,w_0)\to \pi_{n-1}(\breve Q_0,C,w_0)\to \pi_{n-2}(C,w_0)\to
\pi_{n-2}(\breve Q_0,w_0)\cdots  $$

\noindent Here, $\pi_{n-1}(C,w_0)$ need not be trivial, but $C$ bounds a ball in $\breve Q_0$, hence the image of $\pi_{n-1}(C,w_0)$ in
$\pi_{n-1}(\breve Q_0,w_0)$ is trivial.  On the other hand, $\pi_{n-2}(\breve Q_0,w_0)=0$ since $\breve Q_0$ has the
homotopy type of $S^{n-1}$.  Thus we obtain the sequence:

$$\cdots 0\to\integers \to \pi_{n-1}(\breve Q_0,C,w_0)\to \integers\to 0\cdots   $$

\noindent This sequence splits, since we can define a homomorphism\hfil \break$\pi_{n-2}(C,w_0)\to \pi_{n-1}(\breve
Q_0,C,w_0)$ which takes the generator $[C]$ to the class $[E_0]$, and the composition with the boundary operator is
the identity. Thus $\pi_{n-1}(\breve Q_0,C,w_0)\approx \integers\oplus \integers$.  At first glance it appears that the
lift $\tilde f$ of $f$ could take the class $\beta= [\tilde i|_{E_0}]$ to a different class obtained by adding a multiple
of the image of the generator of $\pi_{n-1}(\breve Q_0,w_0)$.  We will show this is impossible.  

If the universal cover $S^{n-1}$ of $V_0$ is a
degree
$s$ finite cover, then we observe that the  generator $\alpha$ of $\pi_{n-1}(\breve Q_0,w_0)$ is
represented by a degree
$s$ covering map $\varphi_1:S^{n-1}\to \tilde i(V_0\times 1)$, whose image includes $s$ copies of $E_0$.  Specifically,
$\varphi_1$ is the covering map $\widetilde{\breve Q}_0\to \breve Q_0$ restricted to the lift to  $\widetilde{\breve Q}_0$ of $V_0\times 1=\tilde
i(V_0\times 1)$.  The map  $\varphi_1$ is homotopic to a similar map $\varphi_0:S^{n-1}\to \tilde i(V_0\times 0)$.  Since $\tilde f$ preserves
orientation on $V_0\times 0$, it takes the homotopy class of $\varphi_0$ to itself, so $\tilde f_*(\alpha)=\alpha$ in $\pi_{n-1}(\breve Q_0)$.  It
follows that $\tilde f$ also preserves the homotopy class of $\varphi_1$ viewed as an element of $\pi_{n-1}(\breve Q_0,C,w_0)$, since
$\tilde f(C)=C$.  We denote the image of $\alpha$ in the relative homotopy group by the same symbol $\alpha$.  The preimage $\varphi_1\inverse(E_0)$
consists of $s$ discs.  Corresponding to the decomposition of $S^{n-1}$ into these $s$ discs and the complement $X$ of their union, there is an equation
in $\pi_{n-1}(\breve Q_0,C,w_0)$: $\alpha=x+s\beta$, where x is represented by $\varphi_1|_X$ and $\beta=[\tilde i|_{E_0}]$.  To interpret $\varphi_1|_X$
as an element of relative
$\pi_{n-1}$, we can use the relative Hurewicz theorem, which says that $\pi_{n-1}(\breve Q_0,C,w_0)\approx H_{n-1}(\tilde Q_0,C)$.

Applying
$\tilde f$ to the equation $\alpha=x+s\beta$, we obtain the equation $\tilde f_*(\alpha)=\tilde f_*(x)+s\gamma$, where 
$\gamma=\tilde f_*(\beta)=[\widetilde
{f\circ i}|_{E_0}]$.  Since ${\tilde f}_*(\alpha)=\alpha$, and also $\tilde f_*(x)=x$ since $\tilde f$ is the identity on $(V_0\times 1)-E_0$, we obtain
the equation $\alpha=x+s\gamma$.  Combining with the earlier equation $\alpha=x+s\beta$, we obtain $s\beta=s\gamma$.  Since this holds in the torsion
free group $\integers\oplus \integers$, it follows that $\beta=\gamma$.  

This completes the proof of the claim.
\end{proof}

From the claim, we now know that 
$\widetilde {f\circ i}|_{E_0}$ is homotopic to $\tilde i|_{E_0}$.  We perform the homotopy rel boundary
on $\widetilde {f\circ i}|_{E}$ and extend to $V_0\times I$ rel $[\bdry (V_0\times I)-E_0]$.  We obtain a map
$\tilde i':V_0\times I\to \breve Q$ homotopic to $\widetilde {f\circ i}$ rel $(\bdry (V_0\times
I)-E_0)$ such that 
$\tilde i=\tilde i'$ on $V_0\times 1$.

Pasting the two maps $\tilde i$ and $\tilde i'$ so that the domain is two copies of $V_0\times I$ doubled on $V_0\times 1$, we obtain a new map
$H:V_0\times I\to \breve Q_0$ with $H|_{V_0\times 0}=\id$ and $H|_{V_0\times 1}=\tilde f|_{V_0\times 0}$.  

Finally
$H$ can be homotoped rel $(V_0\times \bdry I)$ to $ V_0\subset \breve Q_0$ by applying a deformation retraction from
$\breve Q_0$ to $V_0=V_0\times 0$. After this homotopy of $H$, we have $H:V_0\times I\to V_0$ a homotopy in $V_0$ from the identity
to
$f$.  In case $\bdry V_0\ne \emptyset$ in the above argument (for homotoping $H$ to $V_0$), one must perform the homotopy as
a homotopy of pairs, or first homotope $H|_{\bdry V_0\times I}$ to
$\bdry V_0$.  Whenever $\bdry V_0\ne \emptyset$, our construction gives a homotopy of pairs $(f,\bdry f)$ on $(V_0,\bdry V_0)$.

For $Q$ of dimension 4, with components of $V$ irreducible, Grigori Perelman's proof of Thurston's Geometrization Conjecture shows that the universal cover
of each closed component of $V$ is homeomorphic to $\reals^3$ or $S^3$.  Components of $V$ which have non-empty boundary are Haken, possibly
$\bdry$-reducible, and also have contractible universal covers.  By the first part of the lemma, it follows that
$f|_V$ is uniquely determined up to homotopy.  But it is known that automorphisms of such 3-manifolds which are homotopic (homotopic as
pairs) are isotopic, see
\cite{FW:SufficientlyLarge}, \cite{HS:HomotopyIsotopy}, \cite{BO:HomotopyIsotopy}, \cite{BR:HomotopyIsotopy}, \cite{DG:Rigidity},
\cite{GMT:HomotopyHyperbolic}.  It follows that
$f|_V$ is determined by $f|_W$ up to isotopy.
\end{proof}

\section{Proof of Theorem \ref{SequenceThm}.}\label{Sequence}

We begin with a basic lemma which is a consequence of the isotopy extension theorem.

Let $M$ be an $n$-manifold, and let $W$ and $V$ be disjoint compact submanifolds of any (possibly different) dimensions $\le n$.
The submanifolds $W$ and
$V$ need not be properly embedded but they must be submanifolds for which isotopy extension holds, e.g. smooth submanifolds. 
Various mapping class groups of the triple $(M,W,V)$ were defined in the introduction, before the statement of  Lemma \ref{SequenceGeneralLemma}, which
we now prove.

\begin{proof}[Proof of Lemma \ref{SequenceGeneralLemma}]  We need to prove the exactness of the sequence 
$$1\to \H_d(W)\to \H_x(W,V)\to \H_x(V)\to 1,$$

\noindent where the first map is inclusion of an element of $\H_d(W)$ extended to $(W,V)$ by the identity on $V$, the second map is
restriction.  The restriction map takes the group $\H_x(W,V)$ to $\H_x(V)$ and is clearly surjective.  It remains to show that the kernel of this
restriction map is $\H_d(W)$.  If $(f,g)\in \H_x(W,V)$ restricts to the identity on $\H_x(V)$, then $g$ is isotopic to the identity.  By isotopy
extension, we can assume it actually is the identity.  Thus the kernel consists of elements $(f,\id)\in \H_x(W,V)$ such that $(f,\text{id})$ extends
over
$M$, which is the definition of $\H_d(W)$.
\end{proof}

We note that $ \H(W,V)\approx \H(W)\times \H(V)$ is a product, but this does not show that the short exact sequence of Lemma
\ref{SequenceGeneralLemma} splits, and that the group $\H_x(W,V)$ is a product.

\begin{proof} [Proof of Theorem \ref{SequenceThm}]  We will use Lemma \ref{SequenceGeneralLemma}.  We apply the theorem to the triple $(Q,V,W)$,
obtaining the exact sequence $$1\to \H_d(W)\to \H_x(W,V)\to \H_x(V)\to 1.$$  Given an element
$f\in \H_x(W,V)$, this extends to an automorphism of
$(Q,V)$, and by Theorem \ref{FundamentalTheorem},  $f|_V$ is determined up to isotopy, so $H_x(W,V)$ can be replaced by
$H_x(W)$.

To see that $H_x(V)$ in the sequence can be replaced by $\H(V)$, we must show that every element of $\H(V)$ is realized as the restriction of an
automorphism of $(Q,V)$.   We begin by supposing $f$ is an automorphism of $V$.  Observe that $(Q,V)$ can be represented as a (connected)
handlebody
$H$ attached to $V\times 1\subset V\times I$ by 1-handles, with one 1-handle connecting $H$ to each component of $V\times 1$.  In other words,
if $V$ has $c$ components, we obtain $Q$ by identifying $c$ disjointly embedded balls in $\bdry H$ with $c$ balls in $V\times 1$, one in each
component of $V\times 1$.  Denote the set, or union, of balls in $\bdry H$ by $\B_0$ and the set, or union, of balls in $V\times 1$ by $\B_1$.
There is a homeomorphism $\phi:\B_0\to \B_1$ which we use as a glueing map.  Now $(f\times \id)(\B_1)$ is another collection of balls in
$V\times 1$, again with the property that there is exactly one ball in each component of $V\times 1$.  We can then isotope $f\times \id$ to
obtain $F:V\times I\to V\times I$ such that $F(\B_1)=\B_1$.  Then $F$ permutes the balls of $\B_1$, and $\phi\inverse\circ F\circ \phi$ permutes the balls
of $\B_0$.  We use isotopy extension to obtain $G:H\to H$ with $G|_{\B_0}=\phi\inverse\circ F\circ \phi$.  Finally, we can combine $F$ and $G$ by
identifying $F(\B_1)$ with $G(\B_0)$ to obtain an automorphism $g:(Q,V)\to (Q,V)$ with $g|_V=f$.  
\end{proof}

Using ideas in the above proof, we can also prove Theorem \ref{MappingCharacteristicThm}.

\begin{proof}[Proof of Theorem \ref{MappingCharacteristicThm}]
 We apply Lemma \ref{SequenceGeneralLemma}.  Letting $W=\bdry M=\bdry_eQ$ and $V=\bdry_iQ$, where $Q$ is the characteristic compression
body for $M$.  Thus we know that the sequence
$$1\to \H_d(W)\to \H_x(W,V)\to \H_x(V)\to 1$$
\noindent is exact.

In the statement of the theorem, $\H_x(V)$ is defined to be the group of automorphisms of $V$ which extend to $\hat M$, whereas in Lemma
\ref{SequenceGeneralLemma}, it is defined to be the group of automorphisms of $V$ which extend to maps of the triple $(M,W,V)$.  In fact, these are the same,
since any automorphism of $V$ also extends to $Q$, as we showed in the above proof of Theorem \ref{SequenceThm}.

To verify the exactness of the sequence in the statement, it remains to show that $\H_x(W,V)$ is the same as $\H_x(W)$.  To obtain an
isomorphism
$\phi:\H_x(W,V)\to\H_x(W)$, we simply restrict to $W$, $\phi(h)=h|_W$.  Given an automorphism $f$ of $\H_x(W)$ we know there is an
automorphism $g$ of $M$ which extends
$f$.  Since the characteristic compression body is unique up to isotopy, we can, using isotopy extension, replace $g$ by an
isotopic automorphism preserving the pair $(W,V)$, and thus preserving $Q$, then restrict $g$ to $(W,V)$ to obtain an automorphism $\psi(f)$ of
the pair.  This is well-defined by the well-known case of Theorem \ref{FundamentalTheorem}; namely, for an automorphism of a 3-dimensional
compression body, the restriction to the exterior boundary determines the restriction to the interior boundary up to isotopy.  Clearly
$\phi\circ\psi$ is the identity. Now using Theorem
\ref{FundamentalTheorem} again, we can check that $\psi\circ \phi$ is the identity.
\end{proof}

\section{{The Canonical Compression Body\hfill} \newline {for a
3-Manifold.\hfill}
\label{ThreeManifoldSection}}

\begin{defns}  Let $W$ be a compact 3-manifold with no sphere boundary components.  A
{\it symmetric system} $\S$ of disjointly embedded essential spheres is a system with the property
that
$W|\S$ ($W$ cut on $\S$) consists of one ball-with-holes, $\hat V_*$, and other components $\hat V_i$,
$i=1,\ldots k$, bounded by separating spheres $S_i$ of $\S$, each of which is an irreducible 3-manifold $V_i$ with
one open ball removed. 
$\hat V_*$ has one boundary sphere
$S_i$ for each $V_i$, $i=1,\ldots k$, and in addition it has $2\ell$ more spheres, two for each non-separating $S_i$ of $\S$. If
$W$ has
$\ell$
$S^2\times S^1$ prime summands, and $k$ irreducible summands, then there are $\ell$ non-separating $S_i$'s in $\S$
and $k$ separating spheres $S_i$ in $\S$.  The symbol $\S$ denotes the set or the union of the spheres of $\S$.

We denote by $\S'$ the set of duplicates of spheres of $\S$ appearing on $\bdry \hat V_*$. Choose an orientation for $\hat V_*$. Each sphere of $\S'$
corresponds either to  a separating sphere $\S$, which obtains an orientation from $\hat V_*$, or it corresponds to a non-separating sphere of $\S$
with an orientation.

We construct a {\it canonical compression body} associated to the 3-manifold $W$ together with a symmetric system $\S$ by
thickening to replace $W$ by $W\times I$, then attaching a 3-handle along $S_i\times 0$ for each $S_i\in \S$ to obtain a
4-manifold with boundary.  Thus for each 3-handle attached to $S_i\times 0$, we have a 3-ball $E_i$ in the 4-manifold 
with
$\bdry E_i=S_i$.   The boundary of the resulting 4-manifold consists of $W$, one
3-sphere
$V_*$, and the disjoint union of the irreducible non-sphere $V_i$'s.  We cap the 3-sphere $V_*$ with
a 4-ball to obtain the {\it canonical compression body} $Q$ associated to $W$ marked by $\S$.  Note that the
exterior boundary of $Q$ is
$W$ and the interior boundary is $V=\sqcup V_i$.

For each non-separating sphere $S$ in $\S$ we have an {\it associated
separating sphere}
$C(S)$, where $C(S)$ lies in $\hat V_*$ and separates a 3-holed ball containing the duplicates of $S$ coming from cutting on $S$.  
Also, we define $ C(E)$ to be the obvious 3-ball in $Q$ bounded by $C(S)$. 
\end{defns}

We have defined the canonical compression body as an object associated to the 3-manifold $W$ {\it together with} the system $\S$.  This will be our
working definition.  Later, we prove the uniqueness result of Proposition \ref{UniquenessProposition}, which allows us to view the canonical
compression body as being associated to $W$.

In some ways, the dual construction of the canonical compression body may be easier to understand:  Let $V_i$, $i=1\ldots k$ be the
irreducible summands of $W$.  Begin with the disjoint union of the $V_i\times I$ together with a 4-ball $Z$.  For each $i$, identify a 3-ball
$E_i'$ in
$V_i\times 1$ with a 3-ball $E_i''$ in $\bdry Z$ to obtain a ball $E_i$ in the quotient space.  Then, for each $S^2\times S^1$ summand of $W$
identify a 3-ball $E_j'$ in $\bdry Z$ with another 3-ball $E_j''$ in $\bdry Z$ to obtain a disc $E_j$ in the quotient.  (Of course new $E_j'$'s and
$E_j''$'s are chosen to be disjoint in $\bdry Z$ from previous ones.)  The result of all the identifications is the canonical compression body $Q$
for $W$ with balls $E_i$ properly embedded.  We denote by $\E$ the union of $E_i$'s.  The system of spheres $\S=\bdry \E$ is a symmetric system.

We now describe uniqueness for compression bodies associated to a 3-manifold $W$ with a symmetric system.

\begin{proposition}\label{UniquenessProp} If the canonical compression bodies $(Q_1,V_1)$ and
$(Q_2,V_2)$ are associated to a compact 3-manifold
$W$ with symmetric systems $\S_1$ and $\S_2$ respectively, and $f:W\to W$ is any automorphism with $f(\S_1)=\S_2$, then there is a homeomorphism
$g:(Q_1,V_1)\to (Q_2,V_2)$ with $g|_W=f$.  Further, $g|_{V_1}$ is determined up isotopy.

\jerk

\noindent In particular, if $(Q_1,V_1)$ and
$(Q_2,V_2)$ are canonical compression bodies associated to $W$ with the system $\S$, then there is a homeomorphism $g:(Q_1,V_1)\to (Q_2,V_2)$ with
$g|_W=\text{id}$.  Further, $g|_{V_1}$ is determined up isotopy.
\end{proposition}

\begin{proof} Starting with $W\times I$, if we attach 3-handles along the spheres of $\S_1\times 0$, then cap with a 4-ball we obtain $Q_1$.  
Similarly we obtain $Q_2$ from $W\times I$ and $\S_2\times 0$. Starting with $f\times \text{id}:W\times I\to W\times I$, we can then construct the
homeomorphism $g$ by mapping each 3-handle attached on a sphere of $\S_1$ to a 3-handle attached to the image sphere in $\S_2$.  Finally we map the
4-ball in $Q_1$ to the 4-ball in
$Q_2$.  This yields $g$ with $g|_W=f$.  The uniqueness up to isotopy of $g|_{V_1}$ follows from Theorem \ref{FundamentalTheorem}.

Letting $\S_1=\S_2$ and $f=\text{id}$, we get the second statement.
\end{proof}

\begin{lemma}\label{BoundingLemma} The canonical compression body $(Q,V)$ associated to a compact
3-manifold
$W$ with symmetric system $\S$ has the property that every 2-sphere $S\embed W$ bounds a 3-ball properly
embedded in
$Q$.
\end{lemma}

\begin{proof} Make $S$ transverse to $\S$ and consider a curve of intersection innermost in $S$ and
bounding a disc $D$ whose interior is disjoint from $\S$.  If $D$ lies in a holed irreducible $\hat V_i$, $i\ge 1$, then we can
isotope to reduce the number of curves of intersection.  If $D$ lies in $\hat V_*$, then $D$ union a
disc $D'$ in $\S$ forms a sphere bounding a 3-ball in $Q$.  (The 3-ball is embedded in the 4-ball $Z$.)  We replace
$S$ by
$S'=S\cup D'-D$, which can be pushed slightly to have fewer intersections with $\S$ than $S$.  Now if $S'$ bounds a 3-ball,
so does $S$, so by induction on the number of circles of intersection of $S$ with $\S$, we have
proved the lemma if we can start the induction by showing that a 2-sphere $S$ whose intersection with $\S$ is empty must bound a ball 
in $Q$. If such a 2-sphere $S$ lies in $\hat V_*$, then it is obvious that it bounds a 3-ball in $Z$.  If it lies in $\hat V_i$ for some 
$i$, then it must be inessential or isotopic to $S_i$, hence it also bounds a ball (isotopic to $E_i$ if $S$ is isotopic to $S_i$).
\end{proof}

\begin{definition} Let $W$ be a reducible 3-manifold.  Let $S$ be an essential sphere in
$W$ from a symmetric system $\S$.  Or let $S$ be a sphere
associated to a non-separating sphere in $\S$.  An {\it
$\S$-slide automorphism} of
$W$ is an automorphism obtained by cutting on $S$, capping one of the resulting boundary spheres of $W|S$ with a ball $B$ to obtain
$W'$, isotoping the ball along a simple closed curve
$\gamma$ in $W'$ returning the ball to itself, extending the isotopy to $W'$, removing the interior of the ball, and reglueing the
two boundary spheres of
the manifold thus obtained from $W'$.  We emphasize that when $S$ separates a holed irreducible summand, either that irreducible summand can be
replaced by a ball or the remainder of the manifold can be replaced by a ball.  An
$\S$-slide automorphism on a sphere
$S$ can be realized as the boundary of a {\it $\E$-slide automorphism} of
$Q$ rel $V$:  cut $Q$ on the ball $B$ bounded by $S$ ($B$ is either an $E$ in $\E$ or it is a ball $C(E)$ associated to a
non-separating ball $E$ in $\E$) to obtain a compression body
$Q'$  with two spots $B'$ and $B''$, duplicate copies of $B$.  Slide the
duplicate $B'$ in $\bdry_eQ'$ along a closed curve
$\gamma$ in
$\bdry_e Q'$, extending the isotopy to $Q'$.  Then reglue $B'$ to $B''$ to recover the
compression body
$Q$.  

There is another kind of automorphism, called an {\it $\S$-interchanging slide} defined by cutting $W$ on
two separating spheres $S_1$ and $S_2$, either both in $\S$ cutting off holed irreducible summands $\hat V_1$ and $\hat V_2$ from
the remainder $\hat W_0$ of $W$, or both spheres associated to different non-separating $S$ in $\S$ and cutting
off holed $S^2\times S^1$ summands $\hat V_1$ and $\hat V_2$ from $W$.  Cap both boundary spheres in
$\hat W_0$ and slide the capping balls along two paths from one ball to the other, interchanging them.  Then
reglue.  The interchanging slide can also be realized as the boundary of an {\it $\E$-interchanging slide} of $Q$ rel $V$.
We will need to distinguish {\it $\S$-interchanging slides of $S^2\times S^1$ summands}.

We define a further type of automorphism of $W$ or $Q$.  It is called an  $\S$-{\it spin}.  Let $S$ be
a non-separating sphere in $\S$.  We do a half Dehn twist on the separating sphere $C(S)$ to map $S$ to
itself with the opposite orientation.  Again, this spin can be realized as a half Dehn twist on the
3-ball
$C(E)$ bounded by $C(S)$ in $Q$, which is an automorphism of $Q$ rel $V$.  

Finally, we consider Dehn twists on spheres of $\S$.  These are well known to have order at most two in the mapping class group, and we
call them {\it $\S$-Dehn twists}.  These also are boundaries of {\it $\E$-Dehn twists}, automorphisms of $Q$ rel $V$.
\end{definition}

\begin{remark}  Slides, spins, Dehn twists on non-separating spheres, and interchanging slides of $S^2\times S^1$ summands all have
the property that they extend over $Q$ as automorphisms rel $V$, so they represent elements of $\H_d(W)$.  
\end{remark}

We have defined $\S$-slides, and other automorphisms, in terms of the system $\S$ in order to be able to describe the
corresponding $\E$-slides in the canonical maximal compression body $Q$  associated to $W$ with the marking system $\S$.  If we
have not chosen a system $\S$ with the associated canonical compression bodies, there are obvious definitions of {\it slides} (or
interchanging slides, Dehn twists, or spins), without the requirement that the essential spheres involved be in $\S$ or be
associated to spheres in $\S$.

\begin{definition}  Let $\S_1$ and $\S_2$ be two symmetric systems for a 3-manifold $W$ and let $\S_1'$ and $\S_2'$  denote the
corresponding sets of duplicates appearing as boundary spheres of the holed spheres $\hat V_{1*}$ and $\hat V_{2*}$, which are components of $W|\S_1$ and
$W|\S_2$ respectively.  Then an {\it allowable assignment} is a bijection
$a:\S_1'\to \S_2'$ such that
$a$ assigns to a duplicate of a non-separating sphere $S$, another duplicate of a non-separating sphere $a(S)$, and it assigns to the duplicate in $\S_1'$
of a separating sphere of
$\S_1$ a duplicate in $\S_2'$ of a separating sphere of $\S_2$ cutting off a homeomorphic holed irreducible summand.
\end{definition}

If an automorphisms $f:W\to W$ satisfies $f(\S_1)=\S_2$, then $f$ induces an allowable assignment $a:\S_1'\to \S_2'$.

The following is closely related to results due to E. C{\'e}sar de S{\'a} \cite{EC:Automorphisms}, see also
M. Scharlemann in Appendix A of \cite{FB:CompressionBody}, and
D. McCullough in \cite{DM:MappingSurvey} (p. 69).

\begin{lemma}\label{DiscrepantLemma} Suppose $\S$ is a symmetric system for $W$ and $Q$ the
corresponding compression body.  Given any symmetric system
 $\S^\sharp$ for $W$, and an allowable assignment $a:\S'\to {\S^\sharp}'$, there is a composition $f$ of $\S$-slides, $\S$-slide
interchanges, and $\S$-spins such that
$f(\S)=\S^\sharp$, respecting orientation and inducing the map $a:\S'\to {S^\sharp}'$.  
The composition $f$ extends over $Q$ as a
composition of
$\E$-slides,
$\E$-slide interchanges and $\E$-spins, all of which are automorphisms of $Q\rel V$.
\end{lemma}

\begin{proof} $W$ cut on the system $\S$ is a manifold with components $\hat V_i$ and $\hat V_*$; $W$ cut on the system $\S^\sharp$ is a
manifold with components $\hat V_i^\sharp$ and $\hat V_*^\sharp$. Make
$\S^\sharp$ transverse to
$\S$ and consider a curve
$\alpha$ of intersection innermost on
$\S^\sharp$ and bounding a disc $D^\sharp$ in $\S^\sharp$.  The curve $\alpha$ also bounds a disc $D$ in $\S$; in fact, there are two choices
for $D$.  If for either of these choices $D\cup D^\sharp$ bounds a ball, then we can isotope $\S^\sharp$ to eliminate some
curves of intersection.  If $D^\sharp$ is contained in a $\hat V_i$, $i\ge 1$, so that $\hat V_i$ is an irreducible
3-manifold with one hole, then one of the choices of $D$ must give $D\cup D^\sharp$ bounding a ball, so we may assume now
that $D^\sharp\subset \hat V_*$, the holed 3-sphere.  We can also assume that for both choices of $D$, $D\cup D^\sharp$ is a
sphere which (when pushed away from $\bdry \hat V_*$) cuts $\hat V_*$ into two balls-with-holes, neither of which is
a ball.  Now we perform slides to move boundary spheres of $\hat V_*$ contained in a chosen ball-with-holes bounded by $D\cup
D^\sharp$ out of that ball-with-holes.  That is, for each boundary sphere of $\hat V_*$ in the ball-with-holes, perform a slide by
cutting
$W$ on the sphere, capping to replace $\hat V_i$ with a ball, sliding, and reglueing.  One must use a slide path
disjoint from $\S^\sharp$, close to $\S^\sharp$, to return to a different component of $\hat V_*-\S^\sharp$, but this of course does not yield a
closed slide path.  To form a closed slide path, it is necessary to isotope $\S^\sharp\cap V_*$ rel $\bdry V_*$.  Seeing the 
existence of such a slide path requires a little thought; the key is to follow $\S^\sharp$ without intersecting it.
After performing slides to remove all holes of $\hat V_*$ in the holed ball in $V_*$ bounded by $D^\sharp\cup D$,  we can
finally isotope $D^\sharp$ to $D$, and beyond, to eliminate curves of intersection in $\S\cap \S^\sharp$.   Repeating, we can apply
slide automorphisms and isotopies of $\S^\sharp$ to eliminate all curves of intersection.

We know that $\S^\sharp$ also cuts $W$ into manifolds $\hat V_i^\sharp$ homeomorphic to $\hat V_i$, and a holed sphere $\hat V_{*}^\sharp$
and so we conclude that each sphere $S$ of $\S$ bounding $\hat V_i$, which is an irreducible manifold with an open ball removed, is
isotopic to a sphere
$S^\sharp$ in $\S^\sharp$ bounding $\hat V_i^\sharp$ which is also an irreducible manifold with an open ball removed.  Clearly also the
$\hat V_i$ bounded by $S$ is homeomorphic to the $\hat V_i^\sharp$.  Further isotopy therefore makes it possible to
assume $f(\S)=\S^\sharp$.  At this point, however, $f$ does not necessarily assign each $S'\in \S'$ to the desired $a(S')$.  This
can be achieved using interchanging slides and/or spins.  The interchanging slides are used to permute spheres of $\S^\sharp$; the spins are used to
reverse the orientation of a non-separating $S^\sharp \in \S^\sharp$, which switches two spheres of ${\S^\sharp}'$.

Clearly we have constructed $f$ as a composition of automorphisms which extend over $Q$ as $\E$-slides,
$\E$-slide interchanges and $\E$-spins, all of which are automorphisms of $Q\rel V$, hence $f$ also extends as an
automorphism of $Q\rel V$.
\end{proof}

Before proving Proposition \ref{CharacterizationProp}, we mention some known results we will need.  The first result
says that an automorphism of a 3-sphere with $r$ holes which maps each boundary sphere to itself is isotopic to the identity.
The second result concerns an automorphism of an arbitrary 3-manifold:  If it is isotopic to the identity and the isotopy moves the
point along a null-homotopic closed path in the 3-manifold, then the isotopy can be replaced by an isotopy fixing the point at all times.  These 
facts can be proved by applying a result of R. Palais, \cite{RP:RestrictionMaps},\cite{EL:RestrictionMaps}.

\begin{proof}[Proof of Proposition \ref{CharacterizationProp}]  Let $g:W\to W$ be an automorphism.  Then by Lemma
\ref{DiscrepantLemma}, we conclude that there is an an discrepant automorphism $f:(Q,V)\to (Q,V)$ mapping $\S$ to
$g(\S)$ according to the assignment
$\S'\to g(\S)'$ induced by $g$.  Thus $f\inverse\circ g$ is the identity on $\S$, preserving orientation.  Now
we construct $Q$ by attaching 3-handles to $W\times 0\subset W\times I$, a handle attached to each sphere of
$\S$.  The automorphisms
$f\inverse\circ g$ and $f$ clearly both extend over $Q$, hence so does $g$.  This shows $\H(W)=\H_x(W)$.

Showing that $\H_d(W)$ has the indicated generators is more subtle.  Suppose we are given $g:W\to W$ with the property
that
$g$ extends to $Q$ as an automorphism rel $V$.  We also use $g$ to denote the extended automorphism.  We apply the method of proof in
Lemma \ref{DiscrepantLemma} to find an discrepant automorphism $f$, a composition of slides of holed irreducible summands $\hat V_i$,
interchanging slides of $S^2\times S^1$ summands and spins of $S^2\times S^1$ to make
$f\circ g(\S)$ coincide with $\S$, preserving orientation.  Note that $f$ extends over $Q$, and that will be the case below as well, as we
modify
$f$.  Thus we can regard $f$ as an automorphism of $Q$.  Notice also that in the above we do not need interchanging slides of irreducible
summands, since $g$ is the identity on $V$, hence maps each $V_i$ to itself.

At this point $f\circ g$ maps each $\hat V_i$ to itself, and of course it maps $\hat V_*$ to itself, also mapping each boundary sphere of $\hat V_*$ to
itself.  An automorphism of a holed sphere $\hat V_*$ which maps each boundary sphere to itself is isotopic to the identity, as we mentioned above. 
However, we must consider this as an automorphism rel the boundary spheres, and so it may be a composition of twists on boundary spheres.  On
$W-\cup_i
\interior(\hat V_i)$ (a holed connected sum of $S^2\times S^1$'s) this yields some twists on non-separating spheres of $\S$ corresponding to
$S^2\times S^1$ summands.  Precomposing these twists with
$f$, we can now assume that
$f\circ g$ is the identity on
$W-\cup_i \interior(\hat V_i)$. 

It remains to analyze $f\circ g$ restricted to $\hat V_i$.  We do this by considering the restriction of $f\circ g$ to $V_i\times I$, noting that we
should regard this automorphism as an automorphism of the triple $(V_i\times I, E_i', V_i\times 0)$, where $E_i'$ is a duplicate of the disc $E_i$
cutting
$V_i\times I$ from $Q$.   The product gives a homotopy $h_t=f\circ g|_{V_i\times t}$ from $h_0=f\circ g|_{V_i\times 0}=\text{id}$ to $h_1=f\circ
g|_{V_i\times 1}$, where we regard these maps as being defined on $V_i$.  Since $V_i$ is an irreducible 3-manifold, we conclude as in the proof of
Theorem
\ref{FundamentalTheorem}, that $h_0$ is isotopic to $h_1$ via an isotopy $h_t$.  From this isotopy, we obtain a path $\gamma(t)=h_t(p)$ where $p$ is a
point in $E_i'$ when $V_i\times 1$ is identified with $V_i$.  Consider a slide $s$ along $\bar \gamma$:  Cut $Q$ on $E_i$, then slide the
duplicate $E_i'$ along the path $\bar \gamma$, extending the isotopy as usual.  By construction, $s$ is isotopic to the identity via $s_t$,
$s_0=\text{id}$,
$s_1=s$, and $s_t(p)$ is the path $\bar \gamma$.  Combining the isotopies $h_t$ and $s_t$, we see that $(s\circ f\circ g)|_{V_i\times 1}$ is isotopic
to the identity via an isotopy $r_t$ with $r_t(p)$ tracing the path $\gamma\bar\gamma$.  This isotopy can be replaced by an isotopy fixing $p$ or
$E_i'$. Thus after replacing $f$ by $s\circ f$, we have $f\circ g$ isotopic to the identity on $\hat V_i$.  Doing the same thing for every $i$, we
finally obtain $f$ so that $f\circ g$ is isotopic to the identity on $W$.  This shows that
$g$ is a composition of the inverses of the generators used to construct $f$. 
\end{proof}

Now we can give a proof that the canonical compression body is unique in a stronger sense.

\begin{proof} [Proof of Proposition \ref{UniquenessProposition}]  We suppose $(Q_1,V_1)$ is the canonical compression body associated to $W$ with
symmetric system $\S_1$ and $(Q_2,V_2)$ is the canonical compression body associated to $W$ with
symmetric system $\S_2$.  From Proposition \ref{UniquenessProp} we obtain a homeomorphism $h:(Q_1,V_1)\to (Q_2,V_2)$ such that $h(\S_1)=\S_2$.  We are
given a homeomorphism $v:V_1\to V_2$.  From the exact sequence of Theorem \ref{SequenceThm}, we know that there is an automorphism $k:(Q_2,V_2)\to
(Q_2,V_2)$ which extends the automorphism $v\circ (h\inverse|_{V_2}):V_2\to V_2$.  Then $g= k\circ h$ is the required homeomorphism $g:(Q_1,V_1)\to
(Q_2,V_2)$ with the property that $g|_{V_1}=v$.
\end{proof}

\section{Spotted Manifolds.}
\label{SpottedSection}

In this section we present a quick ``compression body point of view" on automorphisms of manifolds with sphere boundary components, or ``spotted
manifolds."

\begin{figure}[ht]
\centering
\psfrag{S}{\fontsize{\figurefontsize}{12}$S$}\psfrag{H}
{\fontsize{\figurefontsize}{12}$H$}\psfrag{a}{\fontsize{\figurefontsize}{12}$\alpha$}
\scalebox{1.0}{\includegraphics{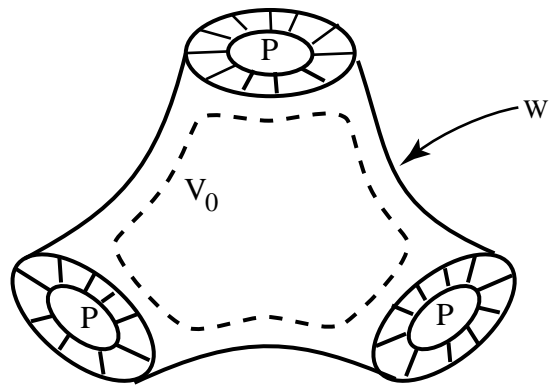}} \caption{\small
 Spotted product.} \label{MapSpotted}
\end{figure}

Suppose $W$ is an $m$-manifold, $m=2$ or $m=3$, with some sphere boundary components.  For such a manifold, there is a canonical collection
$\S$ of essential spheres, namely a collection of spheres isotopic to the boundary spheres.  We can construct the compression body
associated to this collection of spheres, and we obtain an $n=m+1$ dimensional compression body $Q$ whose interior boundary consists
of
$m$-balls, one for each sphere of $\S$, and a closed manifold $V_0$, which is the manifold obtained by capping the sphere boundary
components.  We denote the union of the $m$-balls by $P=\cup P_i$ where each $P_i$ is a ball.  Thus the interior boundary is
$V=V_0\cup P$.  See Figure \ref{MapSpotted}.

Suppose $V_0$ has universal cover which is either $S^m$ or is contractible.  Applying Theorem \ref{SequenceThm}, we obtain the exact
sequence

$$1\to \H_d(W)\to \H_x(W)\to \H(V)\to 1.$$

It remains to interpret the sequence.  Notice that we can view $Q$ as a {\it spotted product} $(V_0\times I,P)$, where $P$ is a
collection of balls in $V_0\times 1$.  However, in our context we have a buffer zone of the form $S^{m-1}\times I\subset A$ 
separating each $P_i$ from $W$. 
We refer to
$W$ as a holed manifold, since it is obtained from the manifold
$V_0$ by removing open balls.  $\H_d(W)$ is the mapping class group of automorphisms which extend to $Q$ as
automorphisms rel
$V=V_0\cup P$.  Either directly, or using Theorem \ref{FundamentalTheorem}, we see that the elements of $\H_d(W)$ are automorphisms $f:W\to
W$ such that the capped automorphism $f_c$ is homotopic, hence isotopic to the identity.  

The mapping class group
$\H_x(W)$ can be identified with $\H(W)$.  This is because any sphere of $\S$ is mapped up to isotopy by any element $f\in \H(W)$ to another sphere of $\S$,
hence $f$ extends to $Q$.  Finally
$\H(V)$ can be identified with
$\H(V_0)\times
\H(P)$, where
$\H(P)$ is clearly just a permutation group.  We obtain the exact sequence 

$$1\to \H_d(W)\to \H(W)\to \H(V_0)\times \H(P)\to 1,$$

\noindent and we have proved Theorem \ref{SpottedThm}.

\bibliographystyle{amsplain}
\bibliography{ReferencesUO3}

\end{document}